\documentclass[10 pt, a4 paper, twoside]{article}
\usepackage{amsmath,amssymb,graphicx}
\usepackage{longtable,mathrsfs}
\usepackage{hyperref}
\usepackage[mathscr]{euscript}
\usepackage{amsmath, amsthm}
\usepackage[mathscr]{euscript}
\newtheorem{thm}{Theorem}

\newtheorem{lem}{Lemma}

\newtheorem{re}{Remark}
\newtheorem{defn}{Definition}
\newcommand{\B}{\mathbb{B}}
\newcommand{\R}{\mathbb{R}}
\newcommand{\N}{\mathbb{N}}
{ \setlength{\oddsidemargin}{-2mm}
  \setlength{\evensidemargin}{6mm} }
  \newlength{\titleright}
\allowdisplaybreaks
\begin{document}
\title{\vspace*{-6pc}{\bf New Approach to Existence of Solution of Weighted Cauchy-type Problem}}

\vspace{1cm}\author{ {Sandeep P Bhairat}\ \\
Department of Mathematics,\\ Institute of Chemical Technology, Mumbai--400 019 (M.S.) India.\\
{Email: sp.bhairat@ictmumbai.edu.in}}

\date{}

\maketitle

{\vspace*{0.5pc} \hrule\hrule\vspace*{2pc}}
\hspace{-0.8cm}{\bf Abstract}\\
We consider a singular fractional differential equation involving generalized Katugampola derivative and obtain the existence and uniqueness of its solution. A scheme for uniformly approximating solution is constructed by using Picard iterative techniques. Illustrative example is also given.\\

\hspace{-0.8cm}{\it \footnotesize {\bf Keywords:}} {\small Fractional integrals and derivatives; Picard iterative technique; singular fractional differential equation.}
\thispagestyle{empty}
\section{Introduction}
In last few decades, the wings of fractional calculus has opened as an emerging trend of applied mathematics with deep applications in almost all branches of science and engineering. At this stage, it ranges to cover the complex problems of real world in physics, control theory, chemical processes and materials, signal and image processing, biological models and dynamical systems. Many researchers devoted to theory and applications of fractional calculus and reported through survey articles \cite{as,ag} and books \cite{ba,hr,kst,fm,pi,skm}.

During the theoretical development of fractional calculus empire, plenty of fractional differential and the corresponding integral operators had come in to existence as well used by timely mathematicians. The Riemann-Liouville, Hadamard, Caputo, Hilfer, Katugampola are the frontiers and their theory became more popular. The investigation of qualitative properties of fractional differential equations is always at the center of development of fractional calculus. The existence and uniqueness of solutions of various fractional differential equations involving these popular operators can be found in \cite{ad,am},\cite{sp3}-\cite{yg},\cite{fj,hh,ke,kh,oo,yy}.

We consider the following weighted Cauchy-type problem
\begin{equation}\label{s3}
\begin{cases}
&\big({^\rho{D}_{a+}^{\alpha,\beta}x} \big)(t)=f(t,x);\qquad t\in\Omega,\,\rho>0,\,0<\alpha<1,\,0\leq\beta\leq1,\\
&\displaystyle\lim_{t\to{a^{+}}}{\big(\frac{t^\rho-a^\rho}{\rho}\big)}^{1-\gamma}x(t)=x_a,\qquad \gamma=\alpha+\beta(1-\alpha),
\end{cases}
\end{equation}
where $f:\Omega\times\R\to\R$ is the given function and ${^\rho{D}_{a+}^{\alpha,\beta}}$ is the generalized Katugampola fractional derivative of order $\alpha$ and type $\beta$, $\Omega=[a,b],0<a<b\le+\infty.$ We prove the existence and uniqueness of solution of Cauchy-type problem \eqref{s3} using equivalent integral equation, properties of gamma function and ratio test (as a convergence criterion). The computable iterative scheme is also constructed to approximate a solution.

The rest of the paper is organised as follows: in next section we list all  definitions and lemmas used throughout the paper. In section 3, we prove equivalent integral equation and existence results followed by illustrative example in section 4. Concluding remarks are given in last section.

\section{Preliminaries}
In this section, we collect some useful definitions and properties from basic fractional calculus \cite{kst,oo,skm}.

As usual $C$ denotes the Banach space of all continuous functions $x:\Omega\to E$ with the superemum (uniform) norm
\begin{equation*}
  {\|x\|}_{\infty}=\sup_{t\in\Omega}{\|x(t)\|}_{E}
\end{equation*}
and $AC(\Omega)$ be the space of absolutely continuous functions from $\Omega$ into $E$. Denote $AC^1(\Omega)-$ the space defined by
\begin{equation*}
AC^1(\Omega)=\bigg\{x:\Omega\to E|\frac{d}{dt}x(t)\in AC(\Omega)\bigg\}.
\end{equation*}
Throughout the paper, let $\delta_{\rho}^n={(t^{\rho-1}\frac{d}{dt})}^{n}, n=[\alpha]+1,$ and mention $[\alpha]$ as integer part of $\alpha.$ Define the space
\begin{equation*}
AC_{\delta_\rho}^{n}=\big\{x:\Omega\to E|{\delta_\rho^{n-1}}x(t)\in AC(\Omega)\big\},\quad n\in\N.
\end{equation*}
Note that $C_{0,\rho}(\Omega)=C(\Omega).$\

Here $L^{p}(a,b), p\geq1,$ is the space of Lebesgue integrable functions on $(a,b).$ The Euler's gamma and beta functions are defined respectively, by
\begin{equation*}
\Gamma(x)=\int_{0}^{+\infty}s^{x-1}e^{-s}ds,\quad \B{(x,y)}=\int_{0}^{1}(1-s)^{x-1}s^{y-1}ds,\quad x>0,y>0.
\end{equation*}
Beta function can be defined through gamma function as: $\B(x,y)=\frac{\Gamma(x)\Gamma(y)}{\Gamma(x+y)},$ for $x>0,y>0,$ \cite{kst}.
\begin{defn}\cite{ki}[Katugampola fractional integral]\label{ki}
Let $\alpha\in{\R}_+,c\in\R$ and $g\in{X_{c}^{p}(a,b)},$ where ${X_{c}^{p}(a,b)}$ is the space of Lebesgue measurable functions. The left-sided Katugampola fractional integral of order $\alpha$ is defined by
\begin{equation*}
({^\rho{I}_{a+}^{\alpha}g})(t)=\int_{a}^{t}s^{\rho-1}{\bigg(\frac{t^\rho-s^\rho}{\rho}\bigg)}^{\alpha-1} \frac{g(s)}{\Gamma(\alpha)}ds,\qquad t>a,\rho>0,
\end{equation*}
where $\Gamma(\cdot)$ is a Euler's gamma function.
\end{defn}
\begin{defn}\cite{kd}
[Katugampola fractional derivative]\label{kd}
Let $\alpha\in{\R}_+\setminus\N$ and $\rho>0.$ The left-sided Katugampola fractional derivative ${^\rho{D}_{a+}^{\alpha}}$ of order $\alpha$ is defined by
\begin{align*}
({^\rho{D}_{a+}^{\alpha}g})(t)&=\delta_\rho^n{({^\rho{I}_{a+}^{n-\alpha}g})(t)}\\
&={\bigg(t^{\rho-1}\frac{d}{dt}\bigg)}^n\int_{a}^{t}s^{\rho-1}{\bigg(\frac{t^\rho-s^\rho}{\rho}\bigg)}^{n-\alpha-1}\frac{g(s)}{\Gamma(n-\alpha)}ds.
\end{align*}
\end{defn}
\begin{defn}\cite{oo}
[Generalized Katugampola fractional derivative]\label{gkd}
The generalized Katugampola fractional derivative of order $\alpha\in(0,1)$ and type $\beta\in[0,1]$ with respect to $t$ and is defined by
\begin{align}\label{gk}
({^\rho{D}_{a+}^{\alpha,\beta}g})(t)=({^\rho{I}_{a+}^{\beta(1-\alpha)}\delta_\rho{^\rho{I}_{a+}^{(1-\beta)(1-\alpha)}g}})(t),\quad\rho>0,
\end{align}
for the function for which right hand side expression exists.
\end{defn}
\begin{re}\label{r1}
The left-sided generalized Katugampola differential operator ${^\rho{D}_{a+}^{\alpha,\beta}}$ can be written as
\begin{equation*}
{^\rho{D}_{a+}^{\alpha,\beta}}={^\rho{I}_{a+}^{\beta(1-\alpha)}}{\delta_\rho}{^\rho{I}_{a+}^{1-\gamma}}={^\rho{I}_{a+}^{\beta(1-\alpha)}}{^\rho{D}_{a+}^{\gamma}},\quad \gamma=\alpha+\beta-\alpha\beta.
\end{equation*}
\end{re}
\begin{re}\cite{oo}
The fractional derivative ${^\rho{D}_{a+}^{\alpha,\beta}}$ is an interpolator of the following fractional derivatives: Hilfer (or genralized R-L) $(\rho\to1)$ \cite{hr}, Hilfer-Hadamard $(\rho\to0^+)$ \cite{hh}, Katugampola $(\beta=0)$ \cite{kd}, Caputo-Katugampola $(\beta=1)$ \cite{am}, Riemann-Liouville $(\beta=0,\rho\to1)$ \cite{kst}, Hadamard $(\beta=0,\rho\to0^+)$ \cite{kh}, Caputo $(\beta=1,\rho\to1)$ \cite{kst}, Caputo-Hadamard $(\beta=1,\rho\to0^+)$ \cite{ad}, Liouville $(\beta=0,\rho\to1,a=-\infty)$ \cite{kst}.
\end{re}
\begin{lem}\label{lm1}\cite{ki}
[Semigroup property]
If $\alpha,\beta>0,1\leq p\leq\infty,0<a<b<\infty$ and $\rho,c\in\R$ for $\rho\geq c.$ Then, for $g\in{X_{c}^{p}(a,b)}$ the following relation hold:
\begin{equation*}
({^\rho{I}_{a+}^{\alpha}}{^\rho{I}_{a+}^{\beta}g})(t)=({^\rho{I}_{a+}^{\alpha+\beta}g})(t).
\end{equation*}
\end{lem}
\begin{lem}\label{lm2}\cite{oo}\label{l2}
Suppose $t>a,~{^\rho{I}_{a+}^{\alpha}}$ and ${^\rho{D}_{a+}^{\alpha}}$ are as in Definition \ref{ki} and Definition \ref{kd}, respectively. Then the following hold:
\begin{align*}
(i)&~~\bigg({^\rho{I}_{a+}^{\alpha}}{\big(\frac{s^\rho-a^\rho}{\rho}\big)}^{\sigma}\bigg)(t)=\frac{\Gamma(\sigma+1)}{\Gamma(\sigma+\alpha+1)}{\big(\frac{t^\rho-a^\rho}{\rho}\big)}^{\sigma+\alpha},\quad\alpha\geq0,\sigma>0,\\
(ii)&~~\text{for}~\sigma=0,~~\bigg({^\rho{I}_{a+}^{\alpha}}{\big(\frac{s^\rho-a^\rho}{\rho}\big)}^{\sigma}\bigg)(t)=\big({^\rho{I}_{a+}^{\alpha}}1\big)(t)=\frac{{\big(\frac{t^\rho-a^\rho}{\rho}\big)}^{\alpha}}{\Gamma(\alpha+1)},\quad\alpha\geq0,\\
(iii)&~~\text{for}~0<\alpha<1,~~\bigg({^\rho{D}_{a+}^{\alpha}}{\big(\frac{s^\rho-a^\rho}{\rho}\big)}^{\alpha-1}\bigg)(t)=0.
\end{align*}
\end{lem}

The following lemma has great importance in the proof of existence results.
\begin{lem}\label{lm3}\cite{pi} Suppose that $x>0.$ Then $\displaystyle\Gamma(x)=\lim_{m\to+\infty}\frac{m^{x}m!}{x(x+1)(x+2)\cdots(x+m)}.$
\end{lem}
We denote $D=[a,a+h], D_{h}=(a,a+h]$, $E=\{x:|x{\big(\frac{t^\rho-a^\rho}{\rho}\big)}^{1-\gamma}-x_a|\leq b\}$ for $h>0,b>0$ and $t\in{D_h}.$ Here we choose $I=(a,a+l]$ and $J=[a,a+l]$ such that
\begin{equation*}
l=\min\bigg{\{h,{\big(\frac{b}{M}\frac{\Gamma(\alpha)}{\B(\alpha,k+1)}\big)}^{\frac{1}{\mu+k}}\bigg\}},\quad \mu=1-\beta(1-\alpha).
\end{equation*}

A function $x(t)$ is said to be a solution of Cauchy-type problem \eqref{s3}, if there exist $l>0$ such that $x\in C^{0}(a,a+l]$ satisfies the differential equation ${^\rho{D}_{a+}^{\alpha,\beta}}x(t)=f(t,x)$ almost everywhere on $I$ alongwith the initial condition $\displaystyle\lim_{t\to{a^{+}}}{{\big(\frac{t^\rho-a^\rho}{\rho}\big)}^{1-\gamma}}x(t)=x_a.$

To prove the existence of solution of Cauchy-type problem \eqref{s3}, let us make following two hypotheses:
\begin{description}
\item[$(H_1)$]  $(t,x)\to f(t,{(\frac{t^\rho-a^\rho}{\rho})}^{\gamma-1}x(t))$ is defined on ${D}_{h}\times E$ satisfies:
\begin{itemize}
\item[(i)] $x\to f(t,{(\frac{t^\rho-a^\rho}{\rho})}^{\gamma-1}x(t))$ is continuous on $E$ for all $t\in{D_{h}}$,\\
 $t\to f(t,{(\frac{t^\rho-a^\rho}{\rho})}^{\gamma-1}x(t))$ is measurable on $D_{h}$ for all $x\in E;$\
\item[(ii)] there exist $k>(\beta(1-\alpha)-1)$ and $M\geq0$ such that $|f(t,{(\frac{t^\rho-a^\rho}{\rho})}^{\gamma-1}x(t))|\leq M{(\frac{t^\rho-a^\rho}{\rho})}^{k}$ holds for all $t\in D_{h}$ and $x\in E,$
\end{itemize}
\item[$(H_2)$] there exists $A>0$ such that $|f(t,{(\frac{t^\rho-a^\rho}{\rho})}^{\gamma-1}x_1(t))-f(t,{(\frac{t^\rho-a^\rho}{\rho})}^{\gamma-1}x_2(t))|$ $\leq A{(\frac{t^\rho-a^\rho}{\rho})}^{k}|x_1-x_2|,$ for all $t\in I$ and $x_1,x_2\in E.$
\end{description}
\section{Existence and uniqueness of solution}
Here, we prove the results for existence and uniqueness of solution of Cauchy-type porblem \eqref{s3}. Also we develop the iterative scheme to approximate the solution and prove its uniqueness.
\begin{lem}\label{lm4}
Suppose that \textbf{$(H_1)$} holds. Then $x:J\to\R$ is a solution of Cauchy-type porblem \eqref{s3} if and only if $x:I\to\R$ is a solution of the integral equation
\begin{equation}\label{a}
x(t)=x_a{\bigg(\frac{t^\rho-a^\rho}{\rho}\bigg)}^{\gamma-1}+\int_{a}^{t}s^{\rho-1}{\bigg(\frac{t^\rho-s^\rho}{\rho}\bigg)}^{\alpha-1}\frac{f(s,x(s))}{\Gamma(\alpha)}ds.
\end{equation}
\end{lem}
\begin{proof} First we suppose that $x:I\to\R$ is a solution of Cauchy-type problem \eqref{s3}. Then, for all $t\in{I},$ we have $|{\big(\frac{t^\rho-a^\rho}{\rho}\big)}^{1-\gamma}x(t)-x_a|\leq b.$  From \textbf{$(H_1)$}, there exists a $k>(\beta(1-\alpha)-1)$ and $M\geq0$ such that
\begin{equation*}
  |f(t,x(t))|=\big|f\big(t,{\big(\frac{t^\rho-a^\rho}{\rho}\big)}^{\gamma-1}{\big(\frac{t^\rho-a^\rho}{\rho}\big)}^{1-\gamma}x(t)\big)\big|\leq M{\big(\frac{t^\rho-a^\rho}{\rho}\big)}^{k},\quad \text{for all}\quad t\in{I}.
\end{equation*}
Then we have,
\begin{align*}
\bigg{|}\int_{a}^{t}{s^{\rho-1}}{\big(\frac{t^\rho-s^\rho}{\rho}\big)}^{\alpha-1}\frac{f(s,x(s))}{\Gamma(\alpha)}ds\bigg{|}&\leq \int_{a}^{t}s^{\rho-1}{\bigg(\frac{t^\rho-s^\rho}{\rho}\bigg)}^{\alpha-1}M\frac{{\big(\frac{s^\rho-a^\rho}{\rho}\big)}^k}{\Gamma(\alpha)}ds\\
&=M{\bigg(\frac{t^\rho-a^\rho}{\rho}\bigg)}^{\alpha+k}\frac{\B(\alpha,k+1)}{\Gamma(\alpha)}.
\end{align*}
Clearly,
\begin{equation*}
\lim_{t\to a+}{\bigg(\frac{t^\rho-a^\rho}{\rho}\bigg)}^{1-\gamma}\int_{a}^{t}{s^{\rho-1}}{\bigg(\frac{t^\rho-s^\rho}{\rho}\bigg)}^{\alpha-1}\frac{f(s,x(s))}{\Gamma(\alpha)}ds=0.
\end{equation*}
It follows that
\begin{equation*}
x(t)=x_a{\bigg(\frac{t^\rho-a^\rho}{\rho}\bigg)}^{\gamma-1}+\int_{a}^{t}{s^{\rho-1}}{\bigg(\frac{t^\rho-s^\rho}{\rho}\bigg)}^{\alpha-1}\frac{f(s,x(s))}{\Gamma(\alpha)}ds,\quad t\in{I}.
\end{equation*}
Since $k>(\beta(1-\alpha)-1),$ then $x\in{C^{0}(I)}$ is a solution of integral equation \eqref{a}.

On the other hand, we can see that $x:I\to\R$ is a solution of integral equation \eqref{a} implies that $x$ is solution of Cauchy-type problem \eqref{s3} defined on $J.$ The proof is complete.
\end{proof}
\begin{re}
	In hypothesis \textbf{$(H_1)$}, if ${(\frac{t^\rho-a^\rho}{\rho})}^{-k}f(t,{(\frac{t^\rho-a^\rho}{\rho})}^{\gamma-1}x(t))$ is continuous on $D\times E,$ one may choose $\displaystyle M=\max_{t\in{D}}{(\frac{t^\rho-a^\rho}{\rho})}^{-k}f(t,{(\frac{t^\rho-a^\rho}{\rho})}^{\gamma-1}x(t))$ continuous on ${D_h\times E}$ for all $x\in E.$
\end{re}
To prove existence and uniqueness of solution of Cauchy-type problem \eqref{s3}, we choose a Picard function sequence as follows:
\begin{equation}\label{pf}
\begin{cases}
  \phi_0(t)&=x_a{\big(\frac{t^\rho-a^\rho}{\rho}\big)}^{\gamma-1},\qquad t\in{I}, \\
  \phi_n(t)&=\phi_0(t)+\displaystyle\int_{a}^{t}{s^{\rho-1}}{\big(\frac{t^\rho-s^\rho}{\rho}\big)}^{\alpha-1}\frac{f(s,\phi_{n-1}(s))}{\Gamma(\alpha)}ds,\quad t\in{I},\quad n=1,2,\cdots.
  \end{cases}
\end{equation}
Now we state the following existence result.
\begin{thm}\label{tm1}
Suppose that \textbf{$(H_1)$} and \textbf{$(H_2)$} hold. Then Cauchy-type problem \eqref{s3} has unique continuous solution $\displaystyle\phi(t)={(\frac{t^\rho-a^\rho}{\rho})}^{\gamma-1}\lim_{n\to\infty}{(\frac{t^\rho-a^\rho}{\rho})}^{1-\gamma}\phi_{n}(t)$ on $I$ with $\phi_0(t)$ and $\phi_n(t)$ given by \eqref{pf}.
\end{thm}
First we prove the continuity of $\phi_{n}$ given by \eqref{pf} as follow:
\begin{lem}\label{lm5}
Suppose \textbf{$(H_1)$} holds. Then $\phi_n$ is continuous on $I$ and satisfies $\big|{\big(\frac{t^\rho-a^\rho}{\rho}\big)}^{1-\gamma}\phi_n(t)-x_a\big|\leq b.$
\end{lem}
\begin{proof} By \textbf{$(H_1)$}, for all $t\in{D_h}$ and $|x{\big(\frac{t^\rho-a^\rho}{\rho}\big)}^{1-\gamma}-x_a|\leq b,$ we have
\begin{equation*}
\bigg|f\bigg(t,{\bigg(\frac{t^\rho-a^\rho}{\rho}\bigg)}^{\gamma-1}x\bigg)\bigg|\leq M{\bigg(\frac{t^\rho-a^\rho}{\rho}\bigg)}^{k}.
\end{equation*}
For $n=1,$ we have
\begin{equation}\label{l1}
  \phi_1(t)=x_a{\bigg(\frac{t^\rho-a^\rho}{\rho}\bigg)}^{\gamma-1}+\int_{a}^{t}{s^{\rho-1}}{\bigg(\frac{t^\rho-s^\rho}{\rho}\bigg)}^{\alpha-1}\frac{f(s,\phi_{0}(s))}{\Gamma(\alpha)}ds.
\end{equation}
Then
\begin{align*}
\bigg{|}\int_{a}^{t}{s^{\rho-1}}{\bigg(\frac{t^\rho-s^\rho}{\rho}\bigg)}^{\alpha-1}\frac{f(s,\phi_{0}(s))}{\Gamma(\alpha)}ds\bigg{|}&\leq \int_{a}^{t}{s^{\rho-1}}{\bigg(\frac{t^\rho-s^\rho}{\rho}\bigg)}^{\alpha-1}M\frac{{\big(\frac{s^\rho-a^\rho}{\rho}\big)}^k}{\Gamma(\alpha)}ds\\
&=M{\bigg(\frac{t^\rho-a^\rho}{\rho}\bigg)}^{\alpha+k}\frac{\B(\alpha,k+1)}{\Gamma(\alpha)}.
\end{align*}
Clearly, $\phi_1\in{C^{0}(I)}$ and from \eqref{l1}, we have
\begin{align}\label{l2}
  \bigg|{\bigg(\frac{t^\rho-a^\rho}{\rho}\bigg)}^{1-\gamma}\phi_1(t)-x_a\bigg|&\leq {\bigg(\frac{t^\rho-a^\rho}{\rho}\bigg)}^{1-\gamma}M{\bigg(\frac{t^\rho-a^\rho}{\rho}\bigg)}^{\alpha+k}\frac{\B(\alpha,k+1)}{\Gamma(\alpha)}\nonumber\\
  &\leq Ml^{\alpha+k+1-\gamma}\frac{\B(\alpha,k+1)}{\Gamma(\alpha)}.
\end{align}
By induction hypothesis, for $n=m,$ suppose that $\phi_m\in{C^{0}(J)}$ and $|{(\frac{t^\rho-a^\rho}{\rho})}^{1-\gamma}\phi_n(t)-x_a|\leq b$ for all $t\in{J}.$ We obtain
\begin{equation}\label{l3}
  \phi_{m+1}(t)=x_a{\bigg(\frac{t^\rho-a^\rho}{\rho}\bigg)}^{\gamma-1}+\int_{a}^{t}{s^{\rho-1}}{\bigg(\frac{t^\rho-s^\rho}{\rho}\bigg)}^{\alpha-1}\frac{f(s,\phi_{m}(s))}{\Gamma(\alpha)}ds.
\end{equation}
From above discussion, we obtain $\phi_{m+1}(t)\in {C^{0}(I)}$ and by \eqref{l3}, we have
\begin{align*}
   \bigg|{\bigg(\frac{t^\rho-a^\rho}{\rho}\bigg)}^{1-\gamma}\phi_{m+1}(t)-x_a\bigg|&\leq {\bigg(\frac{t^\rho-a^\rho}{\rho}\bigg)}^{1-\gamma}\int_{a}^{t}{s^{\rho-1}}{\bigg(\frac{t^\rho-s^\rho}{\rho}\bigg)}^{\alpha-1}M\frac{{\big(\frac{s^\rho-a^\rho}{\rho}\big)}}{\Gamma(\alpha)}ds\\
  &=M{\bigg(\frac{t^\rho-a^\rho}{\rho}\bigg)}^{\alpha+k+1-\gamma}\frac{\B(\alpha,k+1)}{\Gamma(\alpha)} \\
  &\leq Ml^{\alpha+k+1-\gamma}\frac{\B(\alpha,k+1)}{\Gamma(\alpha)}\leq b.
\end{align*}
Thus, the result holds for $n=m+1.$ By using principle of mathematical induction, the result is true for all $n$ and the proof is complete.
\end{proof}
In the following we prove the convergence of sequence $\phi_{n}(t)$.
\begin{thm}\label{tm2}
  Suppose that \textbf{$(H_1)$} and \textbf{$(H_2)$} hold. Then the sequence $\{{(\frac{t^\rho-a^\rho}{\rho})}^{1-\gamma}\phi_n(t)\}$ is uniformly convergent on $J.$
\end{thm}
\begin{proof} For $t\in{J},$ consider the series
\begin{equation*}
{{\bigg(\frac{t^\rho-a^\rho}{\rho}\bigg)}^{1-\gamma}\phi_0(t)}+{{\bigg(\frac{t^\rho-a^\rho}{\rho}\bigg)}^{1-\gamma}[\phi_1(t)-\phi_0(t)]}+\cdots+{{\bigg(\frac{t^\rho-a^\rho}{\rho}\bigg)}^{1-\gamma}[\phi_n(t)-\phi_{n-1}(t)]}+\cdots.
\end{equation*}
Using relation \eqref{l2} in the proof of Lemma \ref{lm5}, we obtain
\begin{equation*}
  {\bigg(\frac{t^\rho-a^\rho}{\rho}\bigg)}^{1-\gamma}|\phi_1(t)-\phi_0(t)|\leq M{\bigg(\frac{t^\rho-a^\rho}{\rho}\bigg)}^{\alpha+k+1-\gamma}\frac{\B(\alpha,k+1)}{\Gamma(\alpha)}, \qquad t\in J.
\end{equation*}
From Lemma \ref{lm5}, we have
\begin{align*}
{\bigg(\frac{t^\rho-a^\rho}{\rho}\bigg)}^{1-\gamma}|\phi_2(t)-\phi_1(t)|=AM\frac{\B(\alpha,k+1)}{\Gamma(\alpha)}\frac{\B(\alpha,\alpha+2k+2-\gamma)}{\Gamma(\alpha)}
{\bigg(\frac{t^\rho-a^\rho}{\rho}\bigg)}^{2(\alpha+k+1-\gamma)}.
\end{align*}
Now suppose for $n=m$
\begin{equation*}
 {\bigg(\frac{t^\rho-a^\rho}{\rho}\bigg)}^{1-\gamma}|\phi_{m+1}(t)-\phi_m(t)|\leq A^{m}M{\bigg(\frac{t^\rho-a^\rho}{\rho}\bigg)}^{(m+1)(\alpha+k+1-\gamma)}\times{P_m},
\end{equation*}
where
\begin{equation}\label{pn}
P_m=\prod_{i=0}^{m}\frac{\B(\alpha,(i+1)k+i(\alpha+1-\gamma)+1)}{\Gamma(\alpha)}.
\end{equation}
We have
\begin{align*}
{\bigg(\frac{t^\rho-a^\rho}{\rho}\bigg)}^{1-\gamma}&\bigg|\phi_{m+2}(t)-\phi_{m+1}(t)\bigg|\\
&\leq{\bigg(\frac{t^\rho-a^\rho}{\rho}\bigg)}^{1-\gamma}\int_{a}^{t}s^{\rho-1}
{{\bigg(\frac{t^\rho-s^\rho}{\rho}\bigg)}^{\alpha-1}}\frac{|f(s,\phi_{m+1}(s))-f(s,\phi_m(s))|}{\Gamma(\alpha)}ds\\
&\leq\frac{{(\frac{t^\rho-a^\rho}{\rho})}^{1-\gamma}}{\Gamma(\alpha)}
\int_{a}^{t}s^{\rho-1}{{\bigg(\frac{t^\rho-s^\rho}{\rho}\bigg)}^{\alpha-1}}
A{\bigg(\frac{s^\rho-a^\rho}{\rho}\bigg)}^{k}\bigg[{\bigg(\frac{s^\rho-a^\rho}{\rho}\bigg)}^{1-\gamma}
|\phi_{m+1}(s)-\phi_m(s)|\bigg]ds\\
&={A^{m+1}M{\bigg(\frac{t^\rho-a^\rho}{\rho}\bigg)}^{(m+2)(\alpha+k+1-\gamma)}}\prod_{i=0}^{m+1}\frac{\B(\alpha,(i+1)k+i(\alpha+1-\gamma)+1)}{\Gamma(\alpha)}.
\end{align*}
This means the result is true for $n=m+1.$ Using the principal of mathematical induction, result is true for all $n.$ i.e.
\begin{equation}\label{l4}
{\bigg(\frac{t^\rho-a^\rho}{\rho}\bigg)}^{1-\gamma}|\phi_{n+2}(t)-\phi_{n+1}(t)|\leq A^{n+1}Ml^{(n+2)(\alpha+k+1-\gamma)}\prod_{i=0}^{n+1}\frac{\B(\alpha,(i+1)k+i(\alpha+1-\gamma)+1)}{\Gamma(\alpha)}.
\end{equation}
Now to prove convergence of sequence $\phi_n$, we consider the series
\begin{equation*}
\sum_{n=1}^{\infty}u_n=\sum_{n=1}^{\infty}MA^{n+1}l^{(n+2)(\alpha+k+1-\gamma)}\prod_{i=0}^{n+1}\frac{\B(\alpha,(i+1)k+i(\alpha+1-\gamma)+1)}{\Gamma(\alpha)}.
\end{equation*}
We obtain
\begin{align*}
\frac{u_{n+1}}{u_n}&=\displaystyle\frac{MA^{n+2}l^{(n+3)(\alpha+k+1-\gamma)}\prod_{i=0}^{n+2}\frac{\B(\alpha,(i+1)k+i(\alpha+1-\gamma)+1)}{\Gamma(\alpha)}}{MA^{n+1}l^{(n+2)(\alpha+k+1-\gamma)}\prod_{i=0}^{n+1}\frac{\B(\alpha,(i+1)k+i(\alpha+1-\gamma)+1)}{\Gamma(\alpha)}}\\
&=Al^{\alpha+k+1-\gamma}\frac{\Gamma((n+3)k+(n+2)(\alpha+1-\gamma)+1)}{\Gamma((n+3)(k+\alpha)+(n+2)(1-\gamma)+1)}.
\end{align*}
Applying Lemma \ref{lm3}, we obtain
\begin{align*}
\frac{u_{n+1}}{u_n}&= Al^{\alpha+k+1-\gamma}\frac{\displaystyle\lim_{m\to\infty}\frac{m^{(n+3)k+(n+2)(\alpha+1-\gamma)+1}m!}{((n+3)k+(n+2)(\alpha+1-\gamma)+1)\cdots((n+3)k+(n+2)(\alpha+1-\gamma)+m+1)}}{\displaystyle\lim_{m\to\infty}\frac{m^{(n+3)(k+\alpha)+(n+2)(1-\gamma)+1}m!}{((n+3)(k+\alpha)+(n+2)(1-\gamma)+1)\cdots((n+3)(k+\alpha)+(n+2)(1-\gamma)+m+1)}}
\end{align*}
\begin{equation*}
=Al^{\alpha+k+1-\gamma}\bigg[\lim_{m\to\infty}m^{-\alpha}\frac{((n+3)(k+\alpha)+(n+2)(1-\gamma)+1)\cdots((n+3)(k+\alpha)+(n+2)(1-\gamma)+m+1)}
{((n+3)k+(n+2)(\alpha+1-\gamma)+1)\cdots((n+3)k+(n+2)(\alpha+1-\gamma)+m+1)}\bigg].
\end{equation*}
We can see that $$\frac{((n+3)(k+\alpha)+(n+2)(1-\gamma)+1)\cdots((n+3)(k+\alpha)+(n+2)(1-\gamma)+m+1)}
{((n+3)k+(n+2)(\alpha+1-\gamma)+1)\cdots((n+3)k+(n+2)(\alpha+1-\gamma)+m+1)}$$ is bounded for all $m,n.$ Then $\displaystyle\lim_{n\to\infty}\frac{u_{n+1}}{u_n}=0.$ Thus $\displaystyle\sum_{n=1}^{\infty}u_n$ is convergent. Hence the series
\begin{equation*}
{{\bigg(\frac{t^\rho-a^\rho}{\rho}\bigg)}^{1-\gamma}\phi_0(t)}+{{\bigg(\frac{t^\rho-a^\rho}{\rho}\bigg)}^{1-\gamma}[\phi_1(t)-\phi_0(t)]}+\cdots+{{\bigg(\frac{t^\rho-a^\rho}{\rho}\bigg)}^{1-\gamma}[\phi_n(t)-\phi_{n-1}(t)]}+\cdots
\end{equation*}
is uniformly convergent. Therefore the sequence $\{{(\frac{t^\rho-a^\rho}{\rho})}^{1-\gamma}\phi_n(t)\}$ is uniformly convergent on $J.$
\end{proof}
\begin{thm}\label{tm3}
Suppose \textbf{$(H_1)$} and \textbf{$(H_2)$} hold. Then $\displaystyle\phi(t)={\bigg(\frac{t^\rho-a^\rho}{\rho}\bigg)}^{\gamma-1}\lim_{n\to\infty}{\bigg(\frac{t^\rho-a^\rho}{\rho}\bigg)}^{1-\gamma}\phi_n(t)$ is unique continuous solution of integral equation \eqref{a} defined on $J.$
\end{thm}
\begin{proof} Since $\displaystyle\phi(t)={\bigg(\frac{t^\rho-a^\rho}{\rho}\bigg)}^{\gamma-1}\lim_{n\to\infty}{\bigg(\frac{t^\rho-a^\rho}{\rho}\bigg)}^{1-\gamma}\phi_n(t)$ on $J,$ and by Lemma \ref{lm5}, we can have\\ ${(\frac{t^\rho-a^\rho}{\rho})}^{1-\gamma}|\phi(t)-x_0|\leq b.$ Then
\begin{align*}
  |f(t,\phi_{n}(t))-f(t,\phi(t))|\leq A{\bigg(\frac{t^\rho-a^\rho}{\rho}\bigg)}^{k}&|\phi_{n}(t)-\phi(t)|,\quad t\in{I},\\
  {\bigg(\frac{t^\rho-a^\rho}{\rho}\bigg)}^{-k}|f(t,\phi_{n}(t))-f(t,\phi(t))|&\leq A|\phi_{n}(t)-\phi(t)|\to0
\end{align*}
uniformly as $n\to+\infty$ on $I.$ Therefore
\begin{align*}
&{\bigg(\frac{t^\rho-a^\rho}{\rho}\bigg)}^{1-\gamma}\phi(t)=\lim_{n\to\infty}\phi_{n}(t)\\
&\hspace{1cm}=x_0+\frac{{(\frac{t^\rho-a^\rho}{\rho})}^{1-\gamma}}{\Gamma(\alpha)}\int_{a}^{t}s^{\rho-1}{\bigg(\frac{t^\rho-s^\rho}{\rho}\bigg)}^{\alpha-1}{\bigg(\frac{s^\rho-a^\rho}{\rho}\bigg)}^{k}\lim_{n\to\infty}\bigg({\bigg(\frac{s^\rho-a^\rho}{\rho}\bigg)}^{-k}f(s,\phi_{n-1}(s))\bigg)ds\\
&\hspace{1cm}=x_0+\frac{{(\frac{t^\rho-a^\rho}{\rho})}^{1-\gamma}}{\Gamma(\alpha)}\int_{a}^{t}s^{\rho-1}{\bigg(\frac{t^\rho-s^\rho}{\rho}\bigg)}^{\alpha-1}f(s,\phi(s))ds.
\end{align*}
Then $\phi$ is a continuous solution of integral equation \eqref{a} defined on $J.$

To prove uniqueness of solution, if possible, suppose that $\psi(t)$ defined on $I$ is also solution of integral equation \eqref{a}. Then ${(\frac{t^\rho-a^\rho}{\rho})}^{1-\gamma}|\psi(t)|\leq b$ for all $t\in{I}$ and
\begin{equation*}
  \psi(t)=x_0{\bigg(\frac{t^\rho-a^\rho}{\rho}\bigg)}^{\gamma-1}+\int_{a}^{t}s^{\rho-1}{\bigg(\frac{t^\rho-s^\rho}{\rho}\bigg)}^{\alpha-1}f(s,\phi(s))ds,\quad t\in{I}.
\end{equation*}
It is sufficient to prove that $\phi(t)\equiv\psi(t)$ on $I.$ From \textbf{$(H_1)$}, there exists a $k>(\beta(1-\alpha)-1)$ and $M\geq0$ such that
\begin{equation*}
  |f(t,\psi(t))|=\bigg|f\bigg(t,{\big(\frac{t^\rho-a^\rho}{\rho}\big)}^{\gamma-1}{\big(\frac{t^\rho-a^\rho}{\rho}\big)}^{1-\gamma}\psi(t)\bigg)\bigg|\leq M{\big(\frac{t^\rho-a^\rho}{\rho}\big)}^{k},
\end{equation*}
for all $t\in{I}.$ Therefore
\begin{align*}
{\bigg(\frac{t^\rho-a^\rho}{\rho}\bigg)}^{1-\gamma}|\phi_{0}(t)-\psi(t)|=&{\bigg(\frac{t^\rho-a^\rho}{\rho}\bigg)}^{1-\gamma}\bigg|\int_{a}^{t}s^{\rho-1}{\bigg(\frac{t^\rho-s^\rho}{\rho}\bigg)}^{\alpha-1}f(s,\psi(s))ds\bigg|\\
&\leq\frac{{(\frac{t^\rho-a^\rho}{\rho})}^{1-\gamma}}{\Gamma(\alpha)}\int_{a}^{t}s^{\rho-1}{\bigg(\frac{t^\rho-s^\rho}{\rho}\bigg)}^{\alpha-1}M{\bigg(\frac{s^\rho-a^\rho}{\rho}\bigg)}^{k}ds\\
&=\frac{M}{\Gamma(\alpha)}{\bigg(\frac{t^\rho-a^\rho}{\rho}\bigg)}^{\alpha+k+1-\gamma}\frac{\B(\alpha,k+1)}{\Gamma(\alpha)}.
\end{align*}
Furthermore, we have
\begin{align*}
{\bigg(\frac{t^\rho-a^\rho}{\rho}\bigg)}^{1-\gamma}|\phi_{1}(t)-\psi(t)|=&\frac{{(\frac{t^\rho-a^\rho}{\rho})}^{1-\gamma}}{\Gamma(\alpha)}\bigg|\int_{a}^{t}s^{\rho-1}{\bigg(\frac{t^\rho-s^\rho}{\rho}\bigg)}^{\alpha-1}[f(s,\phi_0(s))-f(s,\psi(s))]ds\bigg|\\
&\leq AM \frac{\B(\alpha,k+1)}{\Gamma(\alpha)}\frac{\B(\alpha,\alpha+2k+2-\gamma)}{\Gamma(\alpha)}{\bigg(\frac{t^\rho-a^\rho}{\rho}\bigg)}^{2(\alpha+k+1-\gamma)}.
\end{align*}
By the induction hypothesis, we suppose that
\begin{equation*}
{\bigg(\frac{t^\rho-a^\rho}{\rho}\bigg)}^{1-\gamma}|\phi_{n}(t)-\psi(t)|\leq A^{n}M{\bigg(\frac{t^\rho-a^\rho}{\rho}\bigg)}^{(n+1)(\alpha+k+1-\gamma)}\prod_{i=0}^{n}\frac{\B(\alpha,(i+1)k+i(\alpha+1-\gamma)+1)}{\Gamma(\alpha)}.
\end{equation*}
Then
\begin{align*}
{\bigg(\frac{t^\rho-a^\rho}{\rho}\bigg)}^{1-\gamma}|\phi_{n+1}(t)-\psi(t)|\leq&{\bigg(\frac{t^\rho-a^\rho}{\rho}\bigg)}^{1-\gamma}\bigg|\int_{a}^{t}s^{\rho-1}{\bigg(\frac{t^\rho-s^\rho}{\rho}\bigg)}^{\alpha-1}[f(s,\phi_n(s))-f(s,\psi(s))]ds\bigg|\\
\leq& A^{n+1}M{\bigg(\frac{t^\rho-a^\rho}{\rho}\bigg)}^{(n+2)(\alpha+k+1-\gamma)}\prod_{i=0}^{n+1}\frac{\B(\alpha,(i+1)k+i(\alpha+1-\gamma)+1)}{\Gamma(\alpha)}\\
\leq& A^{n+1}Ml^{(n+2)(\alpha+k+1-\gamma)}\prod_{i=0}^{n+1}\frac{\Gamma((i+1)k+i(\alpha+1-\gamma)+1)}{\Gamma((i+1)(\alpha+k)+i(1-\gamma)+1)}.
\end{align*}
By repeating the same arguments used in the proof of Theorem \ref{tm2}, we obtain the series
\begin{equation*}
  \sum_{n=1}^{\infty}A^{n+1}Ml^{(n+2)(\alpha+k+1-\gamma)}\prod_{i=0}^{n+1}\frac{\Gamma((i+1)k+i(\alpha+1-\gamma)+1)}{\Gamma((i+1)(\alpha+k)+i(1-\gamma)+1)}
\end{equation*}
is convergent. Thus
\begin{equation*}
A^{n+1}Ml^{(n+2)(\alpha+k+1-\gamma)}\prod_{i=0}^{n+1}\frac{\Gamma((i+1)k+i(\alpha+1-\gamma)+1)}{\Gamma((i+1)(\alpha+k)+i(1-\gamma)+1)}\to0 ~~\text{as}~~ n\to\infty.
\end{equation*}
We observe that $\displaystyle\lim_{n\to\infty}{\big(\frac{t^\rho-a^\rho}{\rho}\big)}^{1-\gamma}\phi_n(t)={\big(\frac{t^\rho-a^\rho}{\rho}\big)}^{1-\gamma}\psi(t)$ uniformly on $J.$ Thus $\phi(t)\equiv\psi(t)$ on $I.$
\end{proof}
\textbf{Proof of Theorem 3.1:}
\begin{proof} In the light of Lemma \ref{lm4} and from Theorem \ref{tm3}, one can easily deduce that solution $$\phi(t)={\big(\frac{t^\rho-a^\rho}{\rho}\big)}^{\gamma-1}\lim_{n\to\infty}{\big(\frac{t^\rho-a^\rho}{\rho}\big)}^{1-\gamma}\phi_n(t)$$ is unique continuous solution of Cauchy-type problem \eqref{s3} defined on $I$. Thus the proof is ended here.
\end{proof}
\section{An example.}
Will be provided in revised submission.

\section{Conclusion} The existence and uniqueness of solution for a general class of fractional differential equation is obtained using Picard successive approximations. The function $f(t,x)$ considered without assuming the  monotonic property and the iterative scheme is developed for approximating the solution. With the help of well known convergence criteria, the ratio test, the uniform convergence of solution of the considered Cauchy-type problem is established. Our results essentially improves / generalizes the existing results.

\end{document}